\documentclass[12pt,a4paper]{amsart}
\usepackage{amssymb}
\usepackage{mathrsfs}
\headsep=1cm \numberwithin{equation}{section}
\newtheorem{thm}{Theorem}[section]
\newtheorem{cor}[thm]{Corollary}
\newtheorem{lem}[thm]{Lemma}
\newtheorem{prop}[thm]{Proposition}
\theoremstyle{definition}

\theoremstyle{remark}

\numberwithin{equation}{section}

\newcommand\Supp{\operatorname{Supp}}
\newcommand\Ass{\operatorname{Ass}}

\newcommand\Min{\operatorname{Min}}

\newcommand\Rad{\operatorname{Rad}}
\newcommand\Hom{\operatorname{Hom}}
\newcommand\Ext{\operatorname{Ext}}
\newcommand\grade{\operatorname{grade}}
\newcommand\height{\operatorname{height}}

\begin{document}\title[On the Finiteness properties of  local cohomology modules]
{On the finiteness properties of local cohomology modules for regular local rings}
\author [ M. Sedghi, K.  Bahmanpour and  R.  Naghipour]
{Monireh Sedghi, Kamal Bahmanpour and Reza Naghipour$^{\dag}$}
\address{Department of Mathematics, Azarbaijan  Shahid Madani University, 53714-161, Tabriz, Iran.}%
\email {m\_sedghi@tabrizu.ac.ir}  \email {sedghi@azaruniv.ac.ir}%
\address{Department of Mathematics,  Faculty of Mathematical
Sciences, University of Mohaghegh Ardabili, 56199-11367,  Ardabil, Iran.}
 \email{bahmanpour.k@gmail.com}
\address{Department of Mathematics, University of Tabriz, 51666-16471, Tabriz, Iran.}
\email{naghipour@ipm.ir} \email {naghipour@tabrizu.ac.ir}%
\thanks{ 
This research was in part supported by a grant from Azarbaijan Shahid Madani University  (No. 217/d/917).\\
$^{\dag}$Corresponding author: e-mail: {\it naghipour@ipm.ir} (Reza Naghipour)}%
\begin{abstract}
 Let $\frak a$ denote an ideal in a  regular local (Noetherian) ring $R$  and  let $N$  be a finitely generated $R$-module with support in
$V(\frak a)$. The purpose of this paper is to show that all homomorphic images of  the $R$-modules $\Ext^j_R(N, H^i_{\frak a}(R))$ have only finitely many associated primes,
for all $i, j\geq 0$,  whenever  $\dim R \leq4$ or  $\dim R/ \frak a \leq 3$ and $R$ contains a field. In addition, we show that
if  $\dim R=5$ and $R$ contains a field, then the $R$-modules
$\Ext^j_R(N, H^i_{\frak a}(R))$ have only finitely many  associated primes, for all $i, j\geq 0$.
\end{abstract}
\maketitle
\section{Introduction}
In the present paper we continue the study of the finiteness properties of local cohomology modules for regular local rings.
 An interesting problem in commutative algebra is determining when the set of
associated primes of the $i$-th  local cohomology module $H^i_{\frak a}(R)$
of a Noetherian ring $R$ with support in an ideal $\frak a$ of $R$, is finite. This question was
raised by Huneke in \cite{Hu} at the Sundance Conference in 1990.  Examples given by A. Singh \cite{Si} (in the non-local case) and M. Katzman
\cite{Ka} (in the local case)  show there exist local cohomology modules of Noetherian rings with infinitely many associated primes.
However, in recent years there have been several results showing that this conjecture is true in many situations.  The
first result were obtained by Huneke and Sharp. In fact, Huneke and
Sharp \cite{HS} (in the case of positive characteristic) have shown
that, if $R$ is a regular local ring containing a field, then
$H^i_{\frak a}(R)$ has only finitely many associated primes for all
$i\geq0$. Subsequently, G. Lyubeznik in \cite{Ly1} and \cite{Ly2}
showed this result for unramified regular local rings of mixed
characteristic and in characteristic zero.
Further, Lyubeznik posed the following conjecture:\\

{\bf Conjecture.}   {\it If $R$ is a regular ring and $\frak a$ an ideal of $R$, then the local
cohomology modules $H^i_{\frak a}(R)$ have finitely many associated prime ideals for all} $i\geq0$.\\

While this conjecture remains open in this generality, several nice results are now available, see \cite{BN2,He,Mar}. In lower dimensional
cases, Marley in \cite{Mar} showed that the set of associated prime ideals of the local cohomology modules $H^i_{\frak a}(M)$ is finite if $R$ is a local ring of dimension $d$ and $\frak a$ an ideal of $R$, $M$ is a finitely
generated $R$-module, in the following cases: (1) $d\leq3$; (2) $d\leq 4$ and $R$ has an isolated singularity; (3) $d=5$ and $R$ is an unramified
regular local ring and $M$ is torsion-free.

For a survey of recent developments on finiteness properties of
local cohomology modules, see Lyubeznik’s interesting paper \cite{Ly3}.

 The purpose of this paper is to provide some results concerning the set of associated primes  of the local cohomology modules for a regular local ring, that almost results are extensions of Marley's results on local cohomology modules over the strong ring (i.e. the regular local ring containing a field).   Namely, we show that, for a finitely generated module  $N$ over a regular local ring $R$  with support in
$V(\frak a)$, the $R$-modules $\Ext^j_R(N, H^i_{\frak a}(R))$ are
weakly Laskerian, for all $i, j\geq0$,  whenever $\dim R \leq4$ or $\dim R/ \frak a \leq 3$ and $R$ contains a field.
In addition, we show that if $\dim R=5$, then, for all $i, j\geq0$,  the  $R$-module $\Ext^j_R(N, H^i_{\frak a}(R))$  has only finitely many associated primes,  when $R$ contains a field.

We say that an $R$-module $M$ is said to be {\it weakly Laskerian}
if the set of associated primes of any quotient module of $M$ is finite \cite{DM1}.  Our main result in Section 2  is to establish some finiteness results of local cohomology modules for a regular local ring $R$ with respect to an ideal $\frak a$ with $\dim R/ \frak a\leq3$.  More precisely,  we prove the following.\begin{thm}
Let $R$ be a regular local ring containing a field with $\dim R= d\geq 3$, and let $\frak a$ be an ideal of $R$ such that
$\dim R/ \frak a\leq3$.  Then   $\Ext^j_R(R/ \frak a, H^i_{\frak a}(R))$ is weakly Laskerian, for all $i, j\geq 0$.
\end{thm}
The result of Theorem 1.1 is proved in Section 2. Pursuing this point of view further, we  obtain the following consequence of Theorem 1.1, which is an extension of Marley's result in \cite{Mar}.
\begin{cor}
Let $R$ be a regular local ring of dimension $d\leq4$
containing a field and $\frak a$ an ideal of $R$.  Then, for any finitely
generated $R$-module $N$ with $\Supp(N)\subseteq V(\frak a)$, the $R$-module
$\Ext^n_R(N,H^i_{\frak a}(R))$ is weakly Laskerian, for all integers $n,i\geq 0$.
\end{cor}

It will be shown in  Section 3 that the subjects of  Section 2 can be
used to prove a finiteness result  of local cohomology modules for a regular local ring of dimension 5. In fact, we will generalize the main result of Marley for a regular local ring of dimension 5.   More precisely,  we shall show that:

\begin{thm}
Let $R$ be a regular local ring containing a field with $\dim R=5$ and let $\frak a$ be an ideal of $R$.
Then the set $\Ass_R  \Ext^n_R(N, H^i_{\frak a}(R))$ is finite, for each finitely generated $R$-module $N$
with support in $V(\frak a)$ and for all integers $n, i\geq 0$.
\end{thm}

The proof of Theorem 1.3 is given in Section 3.  The following proposition will be one of our main tools for proving Theorem 1.3.
\begin{prop}
Let $R$ be a regular local ring of dimension $d$,  and let $\frak a$ be an ideal of $R$ with $\height  \frak a=1$.  Then
the set $\Ass_R  \Ext^i_R(N,H^1_{\frak a}(R))$ is finite  for each finitely generated $R$-module $N$
with support in $V(\frak a)$ and for all integers $i\geq 0$.
\end{prop}
As a  consequence  of Theorem 1.3 we  derive the following result.
\begin{thm}
Let $R$ be a regular local ring containing a field such that $\dim R\leq5$, and  let $\frak a$ be an ideal of $R$. Then, for all
integers $n,i\geq 0$ and any finitely generated $R$-module $N$ with support in $V(\frak a)$,   the set
 $\Ass_R \Ext^n_R(N, H^i_{\frak a}(R))$ is finite.
\end{thm}
Finally, in this section we will show that, if $M$ is a finitely generated module over a regular local ring $R$ with $\dim R\leq 4$, then
for  any finitely generated $R$-module $N$ with support in $V(\frak a)$, the $R$-module
$\Ext^n_R(N, H^i_{\frak a}(M))$ is weakly Laskerian,  for all integers $i, n\geq 0$.
 In particular the set $\Ass_R H^i_{\frak a}(M)$ is finite, for all integers $i, n\geq 0$.

 Hartshorne \cite{Ha} introduced the notion of a cofinite module, answering in negative a question of Grothendieck \cite[Expos$\acute{e}$  XIII, Conjecture 1.1]{Gr2}. In fact, Grothendieck  asked if the modules  $\Hom_R(R/\frak a,H^{i}_{\frak a}(M))$ always are finitely generated for any ideal $\frak a$ of $R$
and any finitely generated $R$-module $M$. This is the case
when $\frak a=\frak m$, the maximal ideal in a local ring, since the modules $H^{i}_{\frak m}(M)$ are Artinian. Hartshorne defined an
$R$-module $M$ to be $\frak a$-{\it cofinite} if the support of $M$ is contained in $V(\frak a)$ and  $\Ext^i_R (R/ \frak a, M)$ is finitely generated  for all $i\geq 0$.

In \cite{Zo1}, H. Z\"{o}schinger, introduced  an interesting class of
minimax modules, and he has in \cite{Zo1, Zo2} given many equivalent
conditions for a module to be minimax.  An $R$-module $N$ is said to be a minimax module, if there is a finitely generated submodule $L$ of $N$, such that $N/L$ is Artinian. The class of minimax modules thus includes all finitely generated and all Artinian modules. Also, an $R$-module $M$ is called  $\frak a$-{\it cominimax}
if the support of $M$ is contained in $V(\frak a)$ and  $\Ext^i_R (R/ \frak a, M)$ is  minimax  for all $i\geq 0$. The concept of the $\frak a$-cominimax modules is introduced in  \cite{ANV}.

Throughout this paper, $R$ will always be a commutative Noetherian
ring with non-zero identity and $\frak a$ will be an ideal of $R$.
 For an $R$-module $M$, the
$i$-th local cohomology module of $M$ with respect to $\frak a$ is
defined as $$H^i_{\frak a}(M) = \underset{n\geq1} {\varinjlim}\,\,\text{Ext}^i_R(R/ \frak a^n, M).$$
We refer the readers  to \cite{Gr1} and \cite{BS} for more details on local cohomology.

 We shall use $\Min({\frak a})$ to denote the set of all minimal primes
of $\frak a$. For each $R$-module $L$, we denote by
 ${\rm Ass h}_RL$  the set $\{\frak p\in \Ass
_RL:\, \dim R/\frak p= \dim L\}$.  Also, for any ideal $\frak b$ of $R$, we denote
$\{\frak p \in {\rm Spec}\,R:\, \frak p\supseteq \frak b \}$ by
$V(\frak b)$. Finally, for any ideal $\frak c$ of $R$, {\it the
radical of} $\frak c$, denoted by $\Rad(\frak c)$, is defined to
be the set $\{x\in R \,: \, x^n \in \frak c$ for some $n \in
\mathbb{N}\}$. For any unexplained notation and terminology we refer
the readers to \cite{BS} and \cite{Mat}.
\section{Finiteness of local cohomology modules for ideals of small dimension}
 The purpose of this section is to study the finiteness properties of local cohomology modules for a regular local ring $R$ with respect to an ideal $\frak a$ of $R$ with $\dim R/\frak a\leq 3$. The main goal is Theorem 2.7, which
plays a key role in Section 3. This result extends a main result of T. Marley \cite{Mar}.

The following lemmas and proposition will be needed in the proof of Theorem 2.7.
\begin{thm} {\rm(Huneke-Sharp and Lyubeznik).}
Let $(R,\frak m)$ be a regular local ring containing a field. Then
for each ideal $\frak a$ of $R$ and all integers $i\geq 0$, the set of
associated primes of the local cohomology modules $H^i_{\frak a}(R)$ are
finite.
\end{thm}
\proof  See \cite{HS, Ly1, Ly2}. \qed
\begin{lem}
Let $(R,\frak m)$ be a regular local ring
containing a field with  $\dim R=d\geq 3$. Let $\frak p_1,\dots,\frak p_n $  be
 prime ideals of $R$ such that $\dim R/\frak p_i=3$ for all $i=1, \dots, n$.
Then $H^{d}_{\frak a}(R)=0$ and  $\Supp (H^{j}_{\frak a}(R))$ is  finite,  for $j=d-1, d-2$, where $\frak a:=\bigcap_{i=1}^n \frak p_i$.
\end{lem}
\proof It follows easily from Lichtenbaum-Hartshorne and Grothendieck's
vanishing theorems that $H^{d}_{\frak a}(R)=0$, $\Supp (H^{d-1}_{\frak a}(R))\subseteq \{\frak
m\}$ and for each $\frak q \in \Supp (H^{d-2}_{\frak a}(R))$ we have
$\dim R/\frak q\leq 1$.  Therefore $\dim\Supp (H^{d-2}_{\frak a}(R))\leq 1$,  and so $$\Supp (H^{d-2}_{\frak a}(R))\subseteq
{\rm Assh}_R H^{d-2}_{\frak a}(R) \cup \{\frak m\}.$$   Now it
follows from Lemma 2.1  that  $\Supp (H^{d-2}_{\frak a}(R))$ is a finite set, as required.  \qed

\begin{cor}
Let $(R,\frak m)$ and $\frak a$ be as in Lemma {\rm 2.2.} Then  $\Ext^i_R(R/ \frak a, H^{j}_{\frak a}(R))$ is a weakly Laskerian $R$-module for all $i\geq 0$ and for $j=d-1, d-2$.

\proof The result follows easily from Lemma 2.2. \qed
\end{cor}

The next lemma was proved by Melkersson  for $\frak a$-cofiniteness.  The proof given in \cite[Proposition 3.11]{Me3} can be easily carried for weakly Laskerian modules.

\begin{lem}
Let $R$ be a Noetherian ring, $\frak a$  an ideal of $R$, and $M$ an $R$-module such that $\Ext^{i}_{R}(R/ \frak a, M)$ is a weakly Laskerian $R$-module for all $i$.
If $t$  is a non-negative  integer such that the $R$-module $\Ext^{i}_{R}(R/ \frak a, H^{j}_{\frak a}(M))$ is weakly Laskerian,  for all $i$ and all $j\neq t$, then this is the case also when $j=t$.
\end{lem}

\begin{prop}
Let $(R,\frak m)$ be a regular local ring containing a field with  $\dim R=d\geq 3$.  Let $\frak p_1,\dots, \frak p_n $ be  prime ideals of $R$ such that $\dim R/\frak p_i=3$, for all $i=1, \dots, n$, and let $\frak a:=\bigcap_{i=1}^n \frak p_i$.   Then $\Ext^i_R(R/ \frak a, H^{j}_{\frak a}(R))$ is a weakly Laskerian $R$-module for all $i, j \geq 0.$
\end{prop}
\proof If $H^i_{\frak a}(R)\neq 0$,  then it follows from Lemma 2.2 and \cite[Theorems 6.1.2 and  6.2.7]{BS} that $i\in\{d-3,d-2,d-1\}$.
Whence, the assertion follows from Corollary 2.3 and Lemma 2.4. \qed

The next result was proved by Kawasaki  for finitely generated modules.  The proof given in \cite[Lemma 1]{Ka1}
can be easily carried for weakly Laskerian modules.

\begin{lem}
Let $R$ be a Noetherian ring, $T$ an $R$-module,  and $\frak a$ an ideal of $R$. Then the following conditions are equivalent:

{\rm(i)} $\Ext_R^n(R/ \frak a,T)$ is weakly Laskerian for all $n\geq0$.

{\rm(ii)} For any finitely generated $R$-module $N$ with support in $V(\frak a)$, $\Ext_R^n(N,T)$ is weakly Laskerian for all $n\geq0$.
\end{lem}

We now are prepared to prove the main theorem of this section, which shows that when $R$ is a regular local ring
contains a field  and  $\frak a$  an ideal of $R$ such that $\dim R/ \frak a\leq3$, then  all homomorphic images of the
$R$-modules  $\Ext^j_R(R/ \frak a, H^i_{\frak a}(R))$ have  only finitely many associated primes.
\begin{thm}
  Let $(R,\frak m)$ be a regular local ring containing a field with $\dim R= d\geq 3$, and let $\frak a$ be an ideal of $R$ such that
$\dim R/ \frak a\leq3$.  Then  the $R$-module  $\Ext^j_R(R/ \frak a, H^i_{\frak a}(R))$ is weakly Laskerian, for all $i, j\geq 0$.
\end{thm}

\proof In view of  \cite[Corollaries 2.7 and  3.2]{BN3}, we may assume that $\dim R/ \frak a=3$. Then we have $\grade \frak a=d-3$. Hence, by virtue of \cite[Theorems 6.1.2 and 6.2.7]{BS},
 $H^i_{\frak a}(R)=0$,  whenever $i \not \in \{d-3, d-2, d-1\}$. Moreover,  by \cite[Corollary 3.3]{NS}, the
 set $\Supp (H^{d-1}_{\frak a}(R))$ is finite, in view of Proposition 2.5, it is enough for us to show that the
 $R$-module $\Ext^j_R(R/ \frak a, H^{d-3}_{\frak a}(R))$  is weakly Laskerian. To this end, let
\begin{center}
 $X:=\{\frak p \in {\rm Min}(\frak a) \,|\, {\rm height}\, \frak p= d-3\}$\,\,\,\,\,\,\,\, and \,\,\,\,\,\,  $Y:=\{\frak p \in
{\rm Min}(\frak a) \,|\, {\rm height}\,\frak p \geq d-2\}$.
\end{center}
Then $X\neq\emptyset$. Set $\frak b = \bigcap _{\frak p\in X} \frak p$.
 If $Y = \varnothing$, then $\Rad(\frak a)= \frak b$, and so the assertion follows from Proposition 2.5. Therefore we
may assume that $Y \neq \emptyset$. Let $\frak c=\bigcap _{\frak p\in Y} \frak p$.
Then $\Rad(\frak a)=\frak b\cap \frak c$ and ${\rm height}\,\frak c \geq d-2$.
Next, we show that ${\rm height}\,(\frak b+\frak c) \geq d-1$.
Suppose the contrary is true. Then there exists a prime ideal $\frak q$ of $R$ such that $\frak b+\frak c
\subseteq \frak q$ and ${\rm height}\, \frak q = d-2$. Therefore there
exist $\frak p_1 \in X$ and $\frak p_2\in Y$ such that $\frak p_1+
\frak p_2 \subseteq \frak q$. Since ${\rm height}\,\frak p_2 \geq d-2$,
it follows that $\frak q=\frak p_2$, and so $\frak p_1\varsubsetneqq \frak
p_2$, which is a contradiction (note that $\frak p_1, \frak p_2\in {\rm Min}(\frak a)$).

Consequently,  ${\rm grade}\,(\frak b+\frak c) \geq d-1$ and $\grade \frak c\geq d-2$. Therefore,  it
follows from \cite[Theorem 6.2.7]{BS}, that
\begin{center}
$H^{d-3}_{\frak b+\frak c}(R)=0=H^{d-2}_{\frak b+\frak c}(R)$\,\,\,\,\,\, and \,\,\,\,\,\, $H^{d-3}_{\frak c}(R)=0$.
\end{center}
 It now follows from $\Rad(\frak a)=\frak b\cap \frak c$ and the Mayer-Vietoris sequence (see \cite[Theorem 3.2.3]{BS}), that
 $H^{d-3}_{\frak a}(R)\cong H^{d-3}_{\frak b}(R)$, and so by
Proposition 2.5, the $R$-module  $\Ext^j_R(R/\frak b, H^{d-3}_{\frak a}(R))$  is weakly Laskerian, for all $j\geq0$.

On the other hand, since $\dim R/(\frak b +\frak c)\leq 1$, it is easy to see that $$V(\frak b+\frak c)= {\rm Assh}_R\, R/(\frak b +\frak c)\cup\{\frak m\}.$$
Now, as $$\Supp(\Ext^i_R(R/\frak c,H^{d-3}_{\frak b}(R)))\subseteq V(\frak b+\frak c),$$ it follows that
$\Ext^i_R(R/\frak c,H^{d-3}_{\frak b}(R))$ is a weakly Laskerian $R$-module, for all $i\geq0$. Also, as
$${\rm Supp}({\rm Ext}^i_R(R/\frak b+\frak c,H^{d-3}_{\frak b}(R)))\subseteq V(\frak b+\frak c),$$
analogous to the preceding, we see that the $R$-module $${\rm Ext}^i_R(R/\frak b+\frak c,H^{d-3}_{\frak b}(R))$$ is also weakly Laskerian for all $i\geq0$.
Now, the exact sequence
$$0 \longrightarrow R/\Rad(\frak a) \longrightarrow R/\frak b\oplus R/\frak c  \longrightarrow R/\frak b+\frak c \longrightarrow0,$$
induces the long exact sequence
$$\cdots \longrightarrow \Ext^i_R(R/\frak b,H^{d-3}_{\frak a}(R))\oplus\, \Ext^i_R(R/\frak c,H^{d-3}_{\frak a}(R)) \longrightarrow \Ext^i_R(R/\Rad(\frak a), H^{d-3}_{\frak a}(R))$$$$\longrightarrow \Ext^i_R(R/\frak b+\frak c, H^{d-3}_{\frak a}(R))\longrightarrow \cdots,$$
which shows that the $R$-module  $\Ext^i_R(R/\Rad(\frak a), H^{d-3}_{\frak a}(R))$ is weakly Laskerian, for all $i\geq0$, and so it follows from Lemma 2.6 that, the $R$-module
$\Ext^i_R(R/ \frak a, H^{d-3}_{\frak a}(R))$ is weakly Laskerian, for all $i\geq0$, as required. \qed

\begin{cor}
Let $(R,\frak m)$ be a regular local ring of dimension $d\geq 3$
containing a field, and let $\frak a$ be an ideal of $R$ such that $\dim R/ \frak a\leq3$. Then, for any finitely
generated $R$-module $N$ with $\Supp(N)\subseteq V(\frak a)$, the $R$-module $\Ext^n_R(N,H^i_{\frak a}(R))$ is weakly Laskerian, for all integers $n,i\geq 0$.
 \end{cor}

\proof
The result follows from  Theorem 2.7 and Lemma 2.6. \qed

\begin{cor}
Let $(R,\frak m)$ be a regular local ring of dimension $d\leq4$
containing a field, and $\frak a$ an ideal of $R$.  Then, for any finitely
generated $R$-module $N$ with $\Supp(N)\subseteq V(\frak a)$, the $R$-module
$\Ext^n_R(N,H^i_{\frak a}(R))$ is weakly Laskerian, for all integers $n,i\geq 0$.
\end{cor}
\proof Since $R$ is regular local, so $\dim R/ \frak a=4$ if and only if $\frak a=0$. Thus, the assertion
is clear in the case of $\dim R/ \frak a=4$. Moreover, the case $\dim R/ \frak a=3$ follows from
Corollary 2.8. Also, the case $\dim R/ \frak a=2$ follows from \cite[Theorem 3.1]{BN3} and Lemma 2.4. Finally, if
$\dim R/ \frak a\leq1$, then the assertion follows from \cite[Theorem 2.6]{BN3} and \cite[Lemma 1]{Ka1}.\qed

\section{Finiteness of local cohomology modules for regular  local rings of small dimension}

It will be shown in this section that the subjects of the previous section can be
used to prove the finiteness  of local cohomology modules for a regular local ring $R$  with $\dim R=5$. The
main result is Theorem 3.4.  The following proposition will serve to shorten the proof of that theorem.
The following easy lemma will be used in Proposition 3.2.
\begin{lem}
Let $R$ be a Noetherian ring and let $M'\stackrel{f}\longrightarrow M\stackrel{g}\longrightarrow M''$ be an exact sequence of $R$-modules such that $M'$ is weakly Laskerian and $M''$ has only finitely many associated primes. Then $M$ has  only finitely many associated primes.
\end{lem}
\proof The assertion follows from the exact sequence
 $$0\longrightarrow M'/{\rm Im}\,g \longrightarrow M \longrightarrow {\rm Im}\,g \longrightarrow 0,$$
by applying \cite[Theorem 6.3]{Mat}. \qed

\begin{prop}
Let $(R,\frak m)$ be a $d$-dimensional regular local ring,  and $\frak a$ an ideal of $R$ such that $\height  \frak a=1$.  Then
the set $\Ass_R  \Ext^i_R(N,H^1_{\frak a}(R))$ is finite, for each finitely generated $R$-module $N$
with support in $V(\frak a)$ and for all integers $i\geq 0$.
\end{prop}
\proof Let
\begin{center}
 $X:=\{\frak p \in \Min({\frak a}) \,|\, {\rm height}\, \frak p= 1\}$ and $Y:=\{\frak p \in
\Min({\frak a}) \,|\, {\rm height}\,\frak p \geq 2\}$.
\end{center}
Then $X\neq\emptyset$. Let $\frak b = \bigcap _{\frak p\in X} \frak p$.
 Since $R$ is a UFD, it follows from \cite[Exercise 20.3]{Mat} that $\frak b$ is a
principal ideal, and so there is an element $x\in R$ such that $\frak b=Rx$.
Hence, if $Y = \varnothing$, then $\Rad(\frak a)= Rx$, and so it follows from \cite[Theorem 1]{Ka2} that,
the $R$-module $\Ext^i_R(N,H^1_{\frak a}(R))$ is finitely generated. Thus the set $\Ass_R  \Ext^i_R(N,H^1_{\frak a}(R))$ is finite.

Therefore, we may assume that $Y \neq \emptyset$. Then
$\Rad(\frak a)=Rx\cap \frak c$, where $\frak c=\bigcap _{\frak p\in Y} \frak p$. It is easy to see that
${\rm height}\,\frak c \geq 2$ and
${\rm height}\,(\frak b +\frak c) \geq3$. Whence, by using the Mayer-Vietoris sequence it
yields that $H^1_{\frak a}(R)\cong H^1_{\frak b}(R)$. Therefore, by  \cite[Theorem 1]{Ka2},  the $R$-module
$H^1_{\frak a}(R)$ is $\frak b$-cofinite.

On the other hand, according to Artin-Rees lemma, there exists a
positive integer $n$ such that $$\frak b^nN\cap \Gamma_{\frak b}(N)=0=\frak c^nN\cap \Gamma_{\frak c}(N).$$
We claim that $\frak b^nN\cap \frak c^nN= 0$. To this end, suppose that $\frak q\in {\rm Ass}_R(\frak b^nN\cap \frak c^nN)$. Then
$\frak q=(0:_Ry)$, for some $y(\neq0)\in \frak b^nN\cap \frak c^nN$. As $\Supp(N)\subseteq V(\frak a)$,
it follows that $\frak a\subseteq \frak q$, and so $\frak b\cap \frak c\subseteq \frak q$. Thus $\frak b\subseteq \frak q$
or $\frak c\subseteq \frak q$, and so $y\in \Gamma _{\frak b}(N)$ or $y\in \Gamma _{\frak c}(N)$. Furthermore, since
$y\in \frak b^nN\cap \frak c^nN$, it follows that $y=0$,  which is a contradiction.

Now, since $\frak b^nN\cap \frak c^nN= 0$, the exact sequence
$$0\longrightarrow N \longrightarrow N/\frak b^nN\oplus
N/\frak c^nN \longrightarrow N/({\frak b}^n+\frak c^n)N \longrightarrow0,$$induces the
long exact sequence
$$\cdots \longrightarrow \Ext^i_R(N/\frak b^nN, H^{1}_{\frak a}(R))\oplus\, \Ext^i_R(N/\frak c^nN, H^{1}_{\frak a}(R)) \longrightarrow
\Ext^i_R(N, H^{1}_{\frak a}(R))$$$$\longrightarrow \Ext^{i+1}_R(N/({\frak b}^n+\frak c^n)N,H^{1}_{\frak a}(R)) \longrightarrow\cdots.\,\,\,\,\,\,\,\,\,\,\,\,\,\,\,\,\,\,\,\,\,\,\,\,\,\,\,\,\,\,\,\,\,\,\,\,\, (\dag)$$

Since  $H^1_{\frak a}(R)$ is $\frak b$-cofinite and
$$\Supp (N/({\frak b}^n+\frak c^n)N)\subseteq \Supp(N/\frak b^nN)\subseteq V({\frak b}),$$
it follows from \cite[Lemma 1]{Ka1} that the $R$-modules
\begin{center}
${\rm Ext}^i_R(N/\frak b^nN,H^{1}_{\frak a}(R))$  and   ${\rm Ext}^i_R(N/({\frak b}^n+\frak c^n)N,H^{1}_{\frak a}(R))$
\end{center}
are finitely generated for all $i\geq 0$.

 Next, let $T:=N/\frak c^nN$ and we show that ${\rm Ext}^i_R(T,H^{1}_{\frak a}(R))$ is also $\frak b$-cofinite for all $i\geq 0$. To do this, since
${\rm Supp}(\Gamma_{\frak b}(T))\subseteq V({\frak b})$ it follows from \cite[Lemma 1]{Ka1} that the $R$-module ${\rm Ext}^i_R(\Gamma_{\frak b}(T),H^{1}_{\frak a}(R))$ is finitely generated for each $i\geq 0$.  Hence it is enough to show that the  $R$-module  ${\rm Ext}^i_R(T/\Gamma_{\frak b}(T),H^{1}_{\frak a}(R))$ is finitely generated.
As $T/\Gamma_{\frak b}(T)$ is $\frak b$-torsion-free, we therefore make the additional assumption that $T$ is a $\frak b$-torsion-free $R$-module.
Then in view of  \cite[Lemma 2.11]{BS}, $x$ is a non-zerodivisor on $T$, and so the  exact sequence
$$0\longrightarrow T\stackrel{x} \longrightarrow T \longrightarrow T/xT \longrightarrow 0$$
induces the long exact sequence
 $$\cdots \longrightarrow{\rm Ext}^i_R(T/xT,H^{1}_{\frak a}(R))\longrightarrow {\rm Ext}^i_R(T,H^{1}_{\frak a}(R))
\stackrel{x} \longrightarrow {\rm Ext}^i_R(T,H^{1}_{\frak a}(R))$$$$
\longrightarrow {\rm Ext}^{i+1}_R(T/xT,H^{1}_{\frak a}(R)) \longrightarrow
\cdots.\,\,\,\,\,\,\,\,\,\,\,\,\,\,\,\,\,\,\,\,\,\,\,\,\,\,\,\,\,\,\,\,\,\,\,\,\,(\dag\dag)$$

Since ${\rm Supp}(T/xT)\subseteq
V({\frak b})$ it follows from \cite[Lemma 1]{Ka1} that  ${\rm Ext}^i_R(T/xT,H^{1}_{\frak a}(R))$ is a finitely
generated $R$-module, for all $i\geq 0$, (note that $\frak b=Rx$). Consequently, it follows from the
exact sequence $(\dag\dag)$ that the $R$-modules
\begin{center}
$(0:_{{\rm Ext}^i_R(T,H^{1}_{\frak a}(R))}x)$  and  ${\rm Ext}^i_R(T,H^{1}_{\frak a}(R))/x{\rm Ext}^i_R(T,H^{1}_{\frak a}(R))$
\end{center}
are finitely generated and hence $\frak b$-cofinite. Therefore it follows from Melkersson's result \cite[Corollary 3.4]{Me3}
that the $R$-module ${\rm Ext}^i_R(T,H^{1}_{\frak a}(R))$ is $\frak b$-cofinite for all $i\geq 0$.

Now, let $\Omega:=\Ext^i_R(N/\frak b^nN, H^{1}_{\frak a}(R))\oplus\, \Ext^i_R(N/\frak c^nN, H^{1}_{\frak a}(R))$. Then for every finitely
generated submodule $L$ of $\Omega$ the $R$-module $\Omega/L$ is also $\frak b$-cofinite, and so the set $\Ass_R\Omega/L$
is finite.  Now, it follows from the exact sequence $(\dag)$ and Lemma 3.1 that the set
$\Ass_R{\rm Ext}^i_R(N,H^{1}_{\frak a}(R))$ is finite,  as required.\qed\\

The next lemma was proved in \cite{AB} in the case $R$ is local. The proof given in \cite[Lemma 2.5]{AB} can be easily carried over   Noetherian rings, so that we omit the proof.

\begin{lem}
Let $R$ be a Noetherian ring, $x$ an element of $R$ and $\frak a$ an ideal of $R$ such that
$\frak a\subseteq\Rad(Rx)$. Then, for any finitely generated  $R$-module $M$,  the $R$-homomorphism $H^j_{\frak a}(M)\stackrel{x}
\longrightarrow H^j_{\frak a}(M)$ is an isomorphism, for each $j\geq 2$.
\end{lem}

We are now ready to state and prove the main theorem of this section.

\begin{thm}
Let $(R,\frak m)$ be a five-dimensional regular local ring containing a field and $\frak a$ an ideal of $R$.
Then the set $\Ass_R  \Ext^n_R(N, H^i_{\frak a}(R))$ is finite, for each finitely generated $R$-module $N$
with support in $V(\frak a)$ and for all integers $n, i\geq 0$.
\end{thm}

\proof Since $\dim R/ \frak a=5$ if and only if $\frak a=0$,  the assertion
is clear in this case. Hence we consider the case when $\dim R/ \frak a\leq 4$.  If $\dim R/ \frak a\leq 3$, then the result
follows from the proof of Corollary 2.8. Therefore, we may assume that $\dim R/ \frak a= 4$. Then ${\rm height}(\frak a)=1$
and in view of the Lichtenbaum-Hartshorne vanishing theorem $H^5_{\frak a}(R)=0$. Whence $H^i_{\frak a}(R)=0$ whenever
$i\not\in\{1, 2, 3, 4\}$. Also, since by \cite[Corollary 3.3]{NS}, the set $\Supp (H^4_{\frak a}(R))$ is finite, it follows from
$$\Supp (\Ext^n_R(N, H^4_{\frak a}(R)))\subseteq \Supp (H^4_{\frak a}(R))$$
that the set $\Ass_R  \Ext^n_R(N, H^4_{\frak a}(R))$ is also finite. Consequently, in view of Proposition 3.2, we may consider the cases $i=2, 3$.\\

{\bf Case 1.} $i=2$.\\
Suppose that
\begin{center}
 $X:=\{\frak p \in \Min({\frak a}) \,|\, {\rm height}\, \frak p= 1\}$ and $Y:=\{\frak p \in \Min({\frak a})\,|\, {\rm height}\,\frak p \geq 2\}$.
\end{center}
Then  $X\neq\emptyset$ and $Y \neq \emptyset$. Let $\frak b = \bigcap
_{\frak p\in X} \frak p$ and $\frak c=\bigcap _{\frak p\in Y} \frak p$.
 Since $R$ is a UFD, it follows from \cite[Ex. 20.3]{Mat} that $\frak b$ is a
principal ideal, and so there is an element $x\in R$ such that $\frak b=Rx$.
Moreover,  $\Rad(\frak a)=Rx\cap \frak c$, ${\rm height}\,\frak c \geq 2$, and in view of Proposition 3.2, we have
$H^1_{\frak a}(R)\cong H^1_{\frak b}(R)$. Thus, by \cite[Theorem 1]{Ka2},  $H^1_{\frak a}(R)$ is  a $\frak b$-cofinite module.
As $\frak c\not\subseteq \bigcup_{\frak p\in X} \frak p,$  it follows that there is an element $y\in \frak c$ such
that $y\not\in \bigcup _{\frak p\in X} \frak p$. Now,  in view of \cite[Corollary 3.5]{Sc}, there exists the exact sequence
 $$0\longrightarrow H^{1}_{Ry}(H^{1}_{Rx}(R))\longrightarrow H^{2}_{Rx+Ry}(R)
\longrightarrow H^{0}_{Ry}(H^{2}_{Rx}(R))\longrightarrow 0,$$
and so $H^2_{Rx+Ry}(R)\cong H^1_{Ry}(H^1_{Rx}(R))$, (note that $H^{2}_{Rx}(R)=0$).

Using again \cite[Corollary 3.5]{Sc}, to show that there exists the exact sequence
$$0\longrightarrow H^{1}_{Ry}(H^{1}_{\frak a}(R))\longrightarrow H^{2}_{\frak a+Ry}(R)
\longrightarrow H^{0}_{Ry}(H^{2}_{\frak a}(R))\longrightarrow 0,\,\,\,\,\,\,\,\,\,\,\,\,  (\dag)$$
and so it follows from $H^1_{\frak a}(R)\cong H^1_{\frak b}(R)$ that
$$H^{1}_{Ry}(H^{1}_{\frak a}(R))\cong H^{1}_{Ry}(H^{1}_{\frak b}(R))=H^1_{Ry}(H^1_{Rx}(R))\cong H^2_{Rx+Ry}(R).$$

Also, since $xy\in\Rad(\frak a)$, it follows that the $R$-module $H^{2}_{\frak a}(R)$ is $R(yx)$-torsion. Hence
using Lemma 3.3, it is easy to see  that the $R$-module $H^{2}_{\frak a}(R)$ is
$Ry$-torsion. That is  $H^{0}_{Ry}(H^{2}_{\frak a}(R))=H^{2}_{\frak a}(R)$. Consequently, from the exact sequence $(\dag)$  we get the following exact sequence,
$$0\longrightarrow H^2_{Rx+Ry}(R)\longrightarrow H^{2}_{\frak a+Ry}(R) \longrightarrow H^{2}_{\frak a}(R)\longrightarrow 0. \,\,\,\,\,\,\,\,\,\,\,\,  (\dag\dag)$$

Furthermore, in view of Proposition 3.2, there exists a positive integer $n$ such that the sequence
$$0 \longrightarrow N \longrightarrow N/\frak b^nN\oplus N/\frak c^nN \longrightarrow N/({\frak b}^n+\frak c^n)N \longrightarrow0$$
is exact, and so we obtain the long exact sequence
$$\cdots \longrightarrow \Ext^i_R(N/\frak b^nN, H^{2}_{\frak a}(R))\oplus\, \Ext^i_R(N/\frak c^nN, H^{2}_{\frak a}(R)) \longrightarrow
\Ext^i_R(N, H^{2}_{\frak a}(R))$$$$\longrightarrow \Ext^{i+1}_R(N/({\frak b}^n+\frak c^n)N,H^{2}_{\frak a}(R)) \longrightarrow\cdots.\,\,\,\,\,\,\,\,\,\,\,\,\,\,\, (\dag\dag\dag)$$

Since  $$\Supp (N/({\frak b}^n+\frak c^n)N)\subseteq \Supp(N/\frak b^nN)\subseteq V({\frak b}),$$  it follows from Lemma 3.3 that
$$\Ext^j_R(N/\frak b^nN, H^{2}_{\frak a}(R))=0=\Ext^j_R(N/({\frak b}^n+\frak c^n)N, H^{2}_{\frak a}(R)),$$ for all $j\geq 0$. Whence, the exact sequence
$(\dag\dag\dag)$ implies that
$$\Ext^j_R(N, H^{2}_{\frak a}(R))\cong \Ext^j_R(N/\frak c^nN, H^{2}_{\frak a}(R)).$$

 Next, it is easy to see that ${\rm height}(Rx+Ry)=2$, and so   ${\rm height}(\frak a+Ry)=2$. Thus $\dim R/(\frak a+Ry)=3$, and hence
as  $$\Supp (N/\frak c^nN)\subseteq V(\frak c)\subseteq V(\frak a+Ry),$$  it follows from Corollary 2.8 that the
$R$-module $\Ext^j_R(N/\frak c^nN, H^{2}_{\frak a+Ry}(R))$ is weakly Laskerian, for all $j\geq 0$.

Also, as $\grade (Rx+Ry)=2$, it follows from  \cite[Theorem 3.3.1]{BS} and \cite[Proposition 3.11]{Me3} that,  the $R$-module
$H^2_{Rx+Ry}(R)$ is $Rx+Ry$-cofinite.  Now, by modifying the argument of the proof of Proposition 3.2,
one can see that the $R$-module $\Ext^j_R(N/\frak c^nN, H^{2}_{Rx+Ry}(R))$ is $Rx$-cofinite for all $j\geq 0$. In
particular,   $\Ass_R\Ext^j_R(N/\frak c^nN, H^{2}_{Rx+Ry}(R))$ is a finite set,  for all $j\geq 0$.   Moreover, from the exact sequence
$(\dag\dag)$, we deduce the long exact sequence
$$\cdots\longrightarrow \Ext^j_R(N/\frak c^nN, H^2_{\frak a+Ry}(R)) \longrightarrow
\Ext^j_R(N/\frak c^nN, H^{2}_{\frak a}(R))$$$$\longrightarrow \Ext^{j+1} _R(N/\frak c^nN, H^{2}_{Rx+Ry}(R)) \longrightarrow \cdots .$$

Now using Lemma 3.1 and the above long exact sequence induced, it follows from the isomorphism
$$\Ext^j_R(N, H^{2}_{\frak a}(R))\cong \Ext^j_R(N/\frak c^nN, H^{2}_{\frak a}(R)),$$
that the set $\Ass_R\Ext^j_R(N, H^{2}_{\frak a}(R))$ is finite. \\

{\bf Case 2.} $i=3$.\\
Let  $X\neq\emptyset$, $Y\neq\emptyset$  and $y$ be as in the
case 1. Then using the same argument,  it follows from Lemma 3.3 that
the $R$-modules $H^2_{\frak a}(R)$ and $H^3_{\frak a}(R)$ are $Ry$-torsion. Thus
\begin{center}
$H^1_{Ry}(H^2_{\frak a}(R))=0$ and $H^0_{Ry}(H^3_{\frak a}(R))=H^3_{\frak a}(R)$.
\end{center}
 Therefore, using
the exact sequence
$$0\longrightarrow H^{1}_{Ry}(H^{2}_{\frak a}(R))\longrightarrow H^{3}_{\frak a+Ry}(R)
\longrightarrow H^{0}_{Ry}(H^{3 }_{\frak a}(R))\longrightarrow 0,\,\,\,\,\,\,\,\,\,\,\,\,  (\dag)$$
(see \cite[Corollary 3.5]{Sc}), we obtain that $H^3_{\frak a}(R)\cong H^3_{\frak a+Ry}(R)$.

Moreover,  it follows easily from the exact sequence
$$0 \longrightarrow N \longrightarrow N/\frak b^nN\oplus N/\frak c^nN \longrightarrow N/({\frak b}^n+\frak c^n)N \longrightarrow0,$$
that $${\rm Ext}^j_R(N,H^{3}_{\frak a}(R))\cong {\rm Ext}^j_R(N/\frak c^nN,H^{3}_{\frak a}(R)),$$
for all $j\geq 0$. Hence, for all $j\geq 0$ the $R$-module
$\Ext^j_R(N/\frak c^nN, H^{3}_{\frak a+Ry}(R))$ is weakly Laskerian, and so $\Ext^j_R(N/\frak c^nN, H^{3}_{\frak a}(R))$ is also
a weakly Laskerian $R$-module. Therefore $\Ext^j_R(N, H^{3}_{\frak a}(R))$ is a weakly Laskerian $R$-module, and  hence it has finitely many associated primes, as required.\qed

\begin{thm}
Let $(R,\frak m)$ be a regular local ring containing a field such that $\dim R\leq5$, and let $\frak a$ be an ideal of $R$. Then, for all
integers $n,i\geq 0$ and any finitely generated $R$-module $N$ with support in $V(\frak a)$,   the set
 $\Ass_R \Ext^n_R(N, H^i_{\frak a}(R))$ is finite.
\end{thm}
\proof The assertion follows from Corollary 2.8 and Theorem 3.4.\qed\\

 The final theorem of this section shows that if $R$ is a regular local ring with  $\dim R\leq 4$, then the $R$-module $H^i_{\frak a}(M)$ has finitely many associated primes, for all $i\geq 0$ and for any finitely generated $M$ over $R$.  Recall that, an $R$-module $L$ is called an $\frak a$-{\it cominimax module} \cite{ANV} if  $\Supp(L)\subseteq V(\frak a)$ and  $\Ext^i_R (R/ \frak a, L)$ is  minimax,  for all $i\geq 0$.
\begin{thm}
Let $(R,\frak m)$ be a regular local ring such that $\dim R\leq 4$.  Suppose that $\frak a$ is  an ideal of $R$ and $M$  a finitely generated $R$-module.
Then for all integers $i, j\geq 0$, the $R$-module $\Ext^j(R/ \frak a, H^i_{\frak a}(M))$ is  weakly Laskerian.
\end{thm}
\proof If $\dim R\leq 3$, then by virtue of  \cite[Theorem 2.12]{AB},  the $R$-modules $H^i_{\frak a}(M)$ are
$\frak a$-cominimax,  and so the $R$-module $\Ext^j(R/ \frak a, H^i_{\frak a}(M))$ is weakly Laskerian.  Hence we may assume that $\dim R=4$.

Now if  $\dim R/ \frak a=4$,  then $\frak a=0$ and so the result holds. Also,  case $\dim R/ \frak a\leq2$  follows from \cite[Corollary 3.2]{BN3}.
Hence we may assume that $\dim R/ \frak a=3$. Then we have ${\rm heigth}(\frak a)=1$.

On the other hand, if $\dim M\leq 2$ then in view of \cite[Corollary 2.7]{BN3} and \cite[Proposition 5.1]{Me3},  the $R$-module $H^i_{\frak a}(M)$ is $\frak a$-cofinite,
and so  the result is clear.  Also, in the case of $\dim M=3$, the assertion follows from  \cite[Proposition 5.1]{Me3}, \cite[Corollary 3.3]{NS},  the Grothendieck Vanishing Theorem and Lemma 2.4.

Therefore we may assume that $\dim M=4$. Then in view of  \cite[Proposition 5.1]{Me3}, \cite[Corollary 3.3]{NS},  the Grothendieck Vanishing Theorem and Lemma 2.4, it is enough to show that the $R$-module $\Ext^j(R/ \frak a, H^2_{\frak a}(M))$ is  weakly Laskerian, for all $j\geq0$.  To this end,  let
\begin{center}
 $X:=\{\frak p \in \Min({\frak a}) \,|\, {\rm height}\, \frak p= 1\}$ and $Y:=\{\frak p \in \Min({\frak a})  \,|\, {\rm height}\,\frak p \geq 2\}$.
\end{center}
Then $X\neq\emptyset$. Let $\frak b = \bigcap _{\frak p\in X} \frak p$.
 If $Y = \emptyset$ then $\Rad(\frak a)= \frak b$, and so $H^2_{\frak a}(M)=0$, as required.
Therefore we may assume that $Y \neq \emptyset$. Then
$\Rad(\frak a)=\frak b\cap \frak c$, ${\rm height}\,\frak c \geq 2$, and in view of Theorem 2.7,  ${\rm height}\,(\frak b +\frak c) \geq3$, where $\frak c=\bigcap _{\frak p\in Y} \frak p$.

 Now, since $\Rad(\frak a) = \frak b\cap \frak c$, in view of the Mayer-Vietoris sequence  the sequence
$$H^2_{\frak b +\frak c}(M) \longrightarrow H^{2}_{\frak c}(M) \stackrel{f} \longrightarrow H^{2}_{\frak a}(M)\stackrel{g}\longrightarrow H^{3}_{\frak b +\frak c}(M) \,\,\,\,\,\,\,\,\,\,\,\,\,\,\,\, (\dag)$$
is exact.

Since $\dim R/(\frak b +\frak c)\leq 1$ we have $V(\frak b +\frak c)= \Min(\frak b +\frak c)$, and so as
$$\Supp (H^i_{\frak b +\frak c}(M))\subseteq V(\frak b +\frak c),$$  it follows that  the set $\Supp (H^i_{\frak b +\frak c}(M))$ is finite, for $i=2, 3$.
Therefore, it follows from the exact sequence $(\dag)$ that  the $R$-modules ${\rm Ker}\, f$ and ${\rm Im}\, g$ are weakly Laskerian. Thus,
as   $\dim R/\frak c\leq 2 $ it follows from \cite[Corollary 3.3]{BN3} that  the $R$-module $\Ext^i_R(R/\frak c, H^2_{\frak c}(M)/{\rm Ker}\, f)$ is weakly Laskerian
for all $i\geq0$. Now, the exact sequence
$$0\longrightarrow H^2_{\frak c}(M)/{\rm Ker}\, f \longrightarrow H^2_{\frak c}(M) \longrightarrow {\rm Im}\, g \longrightarrow 0,$$
induces the long  exact sequence
$$\Ext^i_R(R/\frak c, H^2_{\frak c}(M)/{\rm Ker}\, f)  \longrightarrow \Ext^i_R(R/\frak c, H^2_{\frak c}(M)) \longrightarrow \Ext^i_R(R/\frak c, {\rm Im}\, g),$$
for all $i\geq0$.  Consequently, in view of Lemma 2.12,  $\Ext^i_R(R/\frak c, H^2_{\frak c}(M))$ is  a weakly Laskerian $R$-module.

Finally, by using the exact sequence
$$0 \longrightarrow R/\Rad(\frak a) \longrightarrow R/\frak b\oplus R/\frak c  \longrightarrow R/\frak b+\frak c\longrightarrow0,$$
we get the  long exact sequence
$$\cdots\longrightarrow  \Ext^i_R(R/\frak b+\frak c,H^{2}_{\frak a}(M)) \longrightarrow \Ext^i_R(R/\frak b\oplus R/\frak c, H^{2}_{\frak a}(M))$$$$\longrightarrow \Ext^i_R(R/\Rad(\frak a), H^{2}_{\frak a}(M))\longrightarrow \Ext^{i+1}_R(R/\frak b+\frak c,H^{2}_{\frak a}(M)) \longrightarrow  \cdots.$$

Now, by applying Lemma 2.10, we obtain that $$\Ext^i_R(R/\frak b+\frak c, H^{2}_{\frak a}(M))=0=\Ext^i_R(R/\frak b, H^{2}_{\frak a}(M)),$$  for each $i\geq 0$, and therefore
 $${\rm Ext}^i_R(R/\Rad(\frak a), H^{2}_{\frak a}(M))\cong {\rm Ext}^i_R(R/\frak c, H^{2}_{\frak a}(M)),$$
for each $i\geq 0$.

Hence, the $R$-module ${\rm Ext}^i_R(R/\Rad(\frak a), H^{2}_{\frak a}(M))$ is weakly Laskerian,  for each $i\geq 0$. Thus,  it follows from  Lemma 2.6  that, the $R$-module $\Ext^i_R(R/ \frak a, H^{2}_{\frak a}(M))$ is  also weakly Laskerian, for
 all $i\geq0$,  and this  completes the proof. \qed

We end this section with a  result which  is a generalization of Corollary 2.8.

\begin{cor}
Let the situation be as in Theorem {\rm3.6}.  Then for  any finitely generated $R$-module $N$ with support in $V(\frak a)$, the $R$-module
$\Ext^j_R(N, H^i_{\frak a}(M))$ is weakly Laskerian,  for all integers $i, j\geq 0$.
 In particular the set $\Ass_R \Ext^j_R(N, H^i_{\frak a}(M))$ is finite, for all integers $i, j\geq 0$.
\end{cor}

\proof
The assertion follows from Theorem 3.6 and Lemma 2.6.\qed

\begin{center}
{\bf Acknowledgments}
\end{center}
The authors are deeply grateful to the referee for a very careful
reading of the manuscript and many valuable suggestions in improving
the quality of the paper. We also would like to thank from the Azarbaijan Shahid Madani University for the financial support (No. 217/d/917).


\begin{thebibliography}{99}
\bibitem{AB}
N. Abazari and K. Bahmanpour, {\it On the finiteness of Bass numbers of local cohomology modules,} J. Alg. Appl. {\bf10} (2011), 783-791.
\bibitem{ANV}
J. A'zami, R. Naghipour and B. Vakili, {\it Finiteness properties of local cohomology modules for $\frak a$-minimax modules},
Proc. Amer. Math. Soc. {\bf137} (2009), 439-448.
\bibitem{BN2}
K. Bahmanpour and R. Naghipour, {\it Associated primes of local cohomology modules and Matlis duality}, J. Algebra. {\bf320} (2008), 2632-2641.
\bibitem{BN3}
K. Bahmanpour and R. Naghipour, {\it Cofiniteness of local cohomology modules for ideals of small dimension}, J. Algebra. {\bf321} (2009), 1997-2011.
\bibitem{BS}
 M.P. Brodmann and R.Y. Sharp, {\it Local cohomology; an algebraic introduction with geometric applications,} Cambridge University Press, Cambridge, 1998.
 \bibitem{DM1}
K. Divaani-Aazar and A. Mafi, {\it Associated primes of local cohomology modules}, Proc. Amer. Math. Soc. {\bf 133} (2005), 655-660.
\bibitem{Gr1}
A. Grothendieck, {\it Local cohomology,} Notes by R. Hartshorne, Lecture Notes in Math. {\bf862}, Springer, New York, 1966.
\bibitem{Gr2}
A. Grothendieck, Cohomologie local des faisceaux coherents et
th$\acute{e}$or$\acute{e}$mes de lefschetz locaux et globaux (SGA2),
North-Holland, Amsterdam, 1968.
\bibitem{Ha}
R. Hartshorne, {\it Affine duality and cofiniteness}, Invent. Math.
{\bf9} (1970), 145-164.
\bibitem{He}
M. Hellus, {\it On the associated primes of a local cohomology module}, J. Algebra. {\bf237} (2001), 406–419.
\bibitem{Hu}
C. Huneke, {\it Problems on local cohomology, Free resolutions in commutative algebra and algebraic geometry},
Res. Notes Math. {\bf2}(1992), 93-108.
\bibitem{HS}
C. Huneke and R.Y. Sharp, {\it Bass numbers of local cohomology module}, Trans. Amer. Math.  Soc. {\bf339} (1993), 765-779.
\bibitem{Ly1}
G. Lyubeznik, {\it Finiteness properties of local cohomology modules (an application of D-modules to commutative algebra)},
Invent. Math. {\bf 113} (1993), 41-55.
\bibitem{Ly2} G. Lyubeznik, {\it Finiteness properties of local cohomology modules for regular local rings of mixed chartacteristic:
the unramified case,} Comm. Algebra, {\bf28} (2000), 5867-5882.
\bibitem{Ly3}
G. Lyubeznik, {\it A partial survey of local cohomology, local cohomology and its applications}, Lectures Notes in Pure and Appl.
Math. {\bf 226} (2002), 121-154.
\bibitem{Ka}
M. Katzman, {\it An example of an infinite set of associated primes of a local cohomology module}, J. Algebra {\bf 252} (2002), 161-166.
\bibitem{Ka1}
K.I. Kawasaki, {\it On the finiteness of Bass numbers of local cohomology modules}, Proc. Amer. Math. Soc. {\bf 124} (1996), 3275-3279.
\bibitem{Ka2}
K.I. Kawasaki, {\it Cofiniteness of local cohomology modules for principal ideals},  Bull. London Math. Soc. {\bf 30} (1998), 241-246.
\bibitem{Mar}
T. Marley, {\it The associated primes of local cohomology modules over rings of small dimension}, Manuscripta Math. {\bf 104} (2001), 519-525.
\bibitem{MV}
T. Marley and J.C. Vassilev, {\it Cofiniteness and associated primes of local cohomology modules}, J. Algebra {\bf 256} (2002), 180-193.
\bibitem{Mat}
H. Matsumura, {\it Commutative ring theory}, Cambridge Univ. Press, Cambridge, UK, 1986.
\bibitem{Me3}
L. Melkersson, {\it Modules cofinite with respect to an ideal}, J. Algebra, {\bf 285} (2005), 649-668.
\bibitem{NS} R. Naghipour and M. Sedghi, {\it A characterization of Cohen-Macaulay modules and local cohomology},
 Arch. Math. {\bf87} (2006), 303-308.
 \bibitem{Sc}
P. Schenzel, {\it Proregular sequences, local cohomology and completion},  Math. Scand. {\bf92} (2003), 161-180.
\bibitem{Si}
A.K. Singh, {\it P-torsion elements in local cohomology modules}, Math. Res. Lett. {\bf7} (2000), 165-176.
\bibitem{Zo1}
H. Z$\ddot{o}$schinger, {\it Minimax modules}, J. Algebra, {\bf102} (1986), 1-32.
\bibitem{Zo2}
H. Z$\ddot{o}$schinger, {\it $\ddot{U}$ber die maximalbedingung
f$\ddot{u}$r radikalvolle untermoduln}, Hokkaido Math. J. {\bf17} (1988), 101-116.


\end{thebibliography}
\end{document}